\documentclass[8pt,a4paper]{article}
\usepackage{amsfonts}
\usepackage{amsmath}
\usepackage{amssymb}
\usepackage{amsthm}
\usepackage[dvips]{graphicx}
\usepackage{latexsym}
\usepackage{mathrsfs}
\usepackage{amscd}
\usepackage{indentfirst}
\usepackage{textcomp}
\usepackage[numbers,square,sort&compress]{natbib} 
\usepackage[colorlinks=true,allcolors=black ]{hyperref}
\allowdisplaybreaks[4]

\newtheorem{Theorem}{Theorem}[section]

\newtheorem{Corollary}[Theorem]{Corollary}
\newtheorem{Lemma}[Theorem]{Lemma}

\theoremstyle{definition}
\newtheorem{Definition}{Definition}[section]
\numberwithin{equation}{section}

\DeclareMathOperator{\Der}{Der}
\DeclareMathOperator{\BDer}{BDer}
\DeclareMathOperator{\IBDer}{IBDer}

\DeclareMathOperator{\Span}{span}
\renewcommand{\d}{\mathrm{d}}
\newcommand{\BF}{\mathbb{F}}

\newcommand{\BN}{\mathbb{N}}
\newcommand{\BZ}{\mathbb{Z}}

\begin{document}
\title{{\bf  Skew-symmetric super-biderivations of the Hamiltonian superalgebra $H(m,n;\underline{t})$}}
\author{\normalsize \bf Da  Xu$^1$\,\,\,\,Xiaoning  Xu$^1$}
\date{{{\small{  1. School of Mathematics, Liaoning University, Shenyang, 110036,
China  }}}}
\maketitle

\begin{abstract}
This paper aims to study the skew-symmetric super-biderivations of the Hamiltonian superalgebra $H(m,n;\underline{t})$. Let $H$ denote the Hamiltonian Lie superalgebra $H(m,n;t)$ over a field of characteristic $p>2$. Utilizing the abelian subalgebra $T_H$ and the weight space decomposition of $H$ with respect to $T_H$, we show the action of a skew-symmetric super-biderivation on the elements of $T_H$ and some specific elements of $H$. Moreover, we prove that all skew-symmetric super-biderivations of $H$ are inner.
\end{abstract}
\renewcommand{\thefootnote}{\fnsymbol{footnote}}
\footnote[0]{ Project Supported by National Natural Science
Foundation of China (No.11501274) and the Science Research Project of Liaoning Provincial Education Department, China (No.L2015203).
\\ Author Email: lnuxxn@163.com (X. Xu)}

\section{Introduction}

Derivations and generalized derivations, as significant subjects in the study of algebras, have been the focus of research over an extended period. In recent years, there has been a growing interest among researchers in biderivations \cite{Benkovic2009,Bresar2018,chang20191,chang20192,chen2016,cheng2017,du2013,Ghosseiri2013,han2016,tang2017,tang20181,tang20182,tang20183,wang2013,wang2011,zhang2006}. Bi-derivations, originating in real linear spaces \cite{Maksa1987}, underwent their initial investigations in prime rings before further developments in the study of biderivations in algebra (see \cite{Bresar1993,Bresar1995,Vukman1990}, for example). Thereafter, Benkovic \cite{Benkovic2009} and Du \cite{du2013} obtained  that under certain conditions a biderivation is a sum of an extremal and an inner biderivation for triangular algebras and generalized matrix algebras. In addition, Han et al. determined all the skew-symmetric biderivations of $\mathcal{W}(a,b)$ and gave the explicit form of each linear commuting map on $\mathcal{W}(a,b)$ (see \cite{han2016}, for details). Afterwards, Chang \cite{chang20191,chang20192} proved that each anti-symmetric biderivation is inner and that all commuting maps are scalar multiplication maps for restricted Cartan-type Lie algebras $W(n;\underline{1})$, $S(n;\underline{1})$ and $H(n;\underline{1})$.

Similar to the biderivations in Lie algebras, the super-biderivations of Lie superalgebras have garnered widespread attention from scholars. In \cite{xia2016}, Xia et al. proved that any super-skewsymmetric super-biderivation of super-Virasoro algebras is inner and described the linear super-commuting maps on super-Virasoro algebras. Furthermore, Cheng \cite{cheng2019} determined all the super-skewsymmetric super-biderivations of the twisted $N = 2$ superconformal algebra, showing that every super-skewsymmetric super-biderivation is inner and that all the linear super-commuting maps are standard. In addition, Chang \cite{chang2021} and Zhao \cite{zhao2020} proved that each skew-symmetric super-biderivation is inner for the generalized Witt Lie superalgebra $W(m,n;\underline{t})$ and the contact Lie superalgebra $K(m,n;\underline{t})$. In \cite{bai2023}, Bai et al. obtained that all the super-symmetric superbiderivations are zero and that all the super-skew-symmetric superbiderivations are inner for simple modular Lie superealgebra of Witt type and special type.

This paper is devoted to studying the skew-symmetric super-biderivations of the Hamiltonian superalgebra $H(m,n;\underline{t})$. And the paper is arranged as follows. In Section 2, we review the basic definitions concerning the Hamiltonian superalgebra $H(m,n;\underline{t})$. In Section 3, we give the definition of super-biderivations of Lie superalgebras and obtain some useful conclusions about the super-biderivations. In Section 4, we use the method of the weight space decomposition of $H(m,n;\underline{t})$ with respect to the abelian subalgebra $T_H$ to prove that all skew-symmetric super-biderivations of $H(m,n;\underline{t})$ are inner.

\section{Preliminaries}

In this section, a brief review of the Hamiltonian superalgebras $H(m,n;\underline{t})$ is given (see \cite{zhang2005}).

Hereafter $\BF$ is an algebraically closed field of characteristic $p>2$ and $\BZ_2 = \{\bar{0}, \bar{1}\}$ is the additive group of order 2. For a vector superspace $V=V_{\bar{0}}\, \oplus\, V_{\bar{1}}$, we write $\d(x)=\alpha$ for the parity of $x\in V_\alpha,\alpha\in \BZ_2$. If $V=\oplus_{i\in \BZ}V_i$ is a $\BZ$-graded vector space, for $x\in V_i,\,i\in \BZ$, $x$ is a $\BZ$-homogeneous element and its $\BZ$-degree $i$ . Throughout this paper, we should mention that once the symbol $\d(x)$  appears in an expression, it implies that $x$ is a $\BZ_2$-homogeneous  element.

Let $\BN_+$ be the set of positive integers and $\BN$ the set of non-negative integers.
Given $m,n\in\BN_+\backslash\{1\}$. For $\alpha=(\alpha_{1},\alpha_{2},\cdots,\alpha_{m})\in \BN^m$, put $|\alpha|=\sum_{i=1}^{m}\alpha_{i}$. For $\beta=(\beta_{1},\beta_{2},\cdots,\beta_{m})\in \BN^m$, we write $\alpha+\beta=(\alpha_{1}+\beta_{1},\alpha_{2}+\beta_{2},\cdots,\alpha_{m}+\beta_{m})$,  $\binom{\alpha}{\beta}=\prod_{i=1}^m\binom{\alpha_i}{\beta_i}$, $\alpha\leq\beta\Longleftrightarrow\alpha_i\leq\beta_i, i=1,2,\cdots,m$. For $\varepsilon _i:=( \delta _{i1},\delta _{i2},\cdots ,\delta _{im})$, where $\delta_{ij}$ is the Kronecker symbol, we abbreviate $x^{(\varepsilon_{i})}$ to $x_{i}$, $i=1, 2, \cdots, m$.  We call $\mathcal{U}(m)$ a $divided\ power\ algebra$ which denotes the  $\BF$-algebra of power series in the variable $x_{1},x_{2},\cdots,x_{m}$. The following formulas hold in $\mathcal{U}(m)$:
  \[x^{(\alpha)}x^{(\beta)}=\binom{\alpha+\beta}{\alpha}x^{(\alpha+\beta)},\ \forall\ \alpha,\beta\in \BN^m.\]

Let $\Lambda(n)$ denote the $Grassmann\ superalgebra$ over $\BF$ in $n$ variables $x_{m+1},x_{m+2},\cdots,x_s$, where $s=m+n$. Denote the tensor product $\mathcal{U}(m)\otimes\Lambda(n)$ by $\Lambda(m,n)$. Then $\Lambda(m,n)$ have  a $\BZ_2$-gradation induced by the trivial $\BZ_2$-gradation of  $\mathcal{U}(m)$ and the natural $\BZ_2$-gradation of $\Lambda(n)$:
  \[\Lambda(m,n)_{\bar{0}}=\mathcal{U}(m)\otimes\Lambda(n)_{\bar{0}},\qquad \Lambda(m,n)_{\bar{1}}=\mathcal{U}(m)\otimes\Lambda(n)_{\bar{1}}.\]

Obviously, $\Lambda(m,n)$ is an associative superalgebra. For $g\in\mathcal{U}(m), f\in\Lambda(n)$, we simply write $g\otimes f$ as $gf$. The following formulas hold in $\Lambda(m,n)$:
  \[x_ix_j=-x_jx_i, \qquad i,j=m+1,\cdots,s.\]
  \[x^{(\alpha)}x_j=x_jx^{(\alpha)}, \quad\ \forall\ \alpha\in \BN^m,j=m+1,\cdots,s.\]
For $k = 1,\cdots,n$, set
  \[B_k:=\{\langle i_1,i_2,\cdots,i_k\rangle\mid m+1\leq i_1<i_2<\cdots<i_k\leq s\}\]
and $B(n):=\bigcup_{k=0}^nB_k$, where $B_0=\varnothing$. For $u=\langle i_1,i_2,\cdots,i_k\rangle\in B_k$, set $|u|=k, \{u\}=\{i_1,i_2,\cdots i_k\}$ and $x^u=x_{i_1}x_{i_2}\cdots x_{i_k}$. Specially, let $|\varnothing|=0$, $x^\varnothing=1$. It is obvious that $\{x^{(\alpha)}x^u\mid \alpha\in \BN^m,u\in B(n)\}$ is an $\BF$-basis of $\Lambda(m,n)$.

Let $Y_0=\{1,2,\cdots,m\}$, $Y_1=\{m+1,\cdots,s\}$, and $Y=Y_0\cup Y_1$. Let $D_1,D_2,\cdots,D_s$ be the linear transformations of $\Lambda(m,n)$ such that
  \begin{align*}
  	D_i(x^{(\alpha)}x^u)=
  	\begin{cases}
  		x^{(\alpha-\varepsilon_i)}x^u,\quad&\forall\ i\in Y_0,\\
  		x^{(\alpha)}\partial_i(x^u),\quad &\forall\ i\in Y_1,
  	\end{cases}
  \end{align*}
where $\partial_i$ is the special derivation of $\Lambda(n)$. Then $D_1,D_2,\cdots,D_s$ are superderivations of the superalgebra $\Lambda(m,n)$, and $\d(D_i)=\tau(i)$, where
  \begin{align*}
  	\tau(i)=
  	\begin{cases}
  		\bar{0},\qquad \forall\ i\in Y_0,\\
  		\bar{1},\qquad \forall\ i\in Y_1.
  	\end{cases}
  \end{align*}
Let
  \[W(m,n):=\bigg\{\sum_{i=1}^sf_iD_i\mid f_i\in \Lambda(m,n),\ \forall\ i\in Y\bigg\}.\]
Then $W(m,n)$ is an infinite-dimensional Lie superalgebra which is contained in $\Der(\Lambda(m,n))$ and the following formula holds:
  \[ [fD_i,gD_j]=fD_i(g)D_j-(-1)^{\d(fD_i)\d(gD_j)}gD_j(f)D_i, \]
for all $f,g\in\Lambda(m,n)$ and $i,j\in Y$.

Fix two $m$-tuples of positive integers $\underline{t}=(t_1,t_2,\cdots,t_m)$ and $\pi=(\pi_1,\pi_2,\cdots,\pi_m)$,
where $\pi_i=p^{t_i}-1$ for all $i\in Y_0$ and $p$ is the characteristic of the basic field $\BF$. Let
  \[\Lambda(m,n;\underline{t}):= \Span_{\BF}\{x^{(\alpha)}x^u\mid \alpha\in A(m,\underline{t}),u\in B(n)\}.\]
where $A(m,\underline{t} ) =\{\alpha=(\alpha _1,\alpha _2,\cdots ,\alpha _m)\in \BN^m\mid 0\leq \alpha _i\leq \pi _i,i\in Y_0 \}$. Then $\Lambda(m,n;\underline{t})$ is a subalgebra of $\Lambda(m,n)$. Let 
  \[W(m,n;\underline{t}):=\bigg\{\sum_{i=1}^{s}f_iD_i\mid f_i\in \Lambda(m,n;\underline{t}),\ \forall\ i\in Y\bigg\}.\]
Then $W(m,n;\underline{t})$ is a finite-dimensional simple Lie superalgebra, which is called the generalized Witt Lie superalgebra. Note that $W(m,n;\underline{t})$ possesses a $\BZ$-graded structure:
  \[W(m,n;\underline{t})=\bigoplus\limits_{r=-1}^{\xi-1}W(m,n;\underline{t})_r,\]
by letting $W(m,n;\underline{t})_r:=\Span_{\BF}\{x^{(\alpha)}x^uD_j\mid|\alpha|+|u|=r+1,j\in Y\}$ and $\xi:=|\pi|+n$.

Hereafter, suppose $m=2k$ is even and $n=2t$ is even.  For $i\in Y$, put
\begin{align*}
	i'&=\begin{cases}
		i+k,&\quad   1\leq i\leq k, \\
		i-k,& \quad  k<i\leq 2k, \\
		i+t,& \quad  2k< i\leq 2k+t, \\
		i-t,&  \quad  2k+t< i\leq s.
	\end{cases}\\
	\sigma(i) &=\begin{cases}
		1,&\quad  1\leq i\leq k, \\
		-1,&\quad  k<i\leq 2k,\\
		1,& \quad  2k< i\leq s.
	\end{cases}
\end{align*}
Define a linear mapping $D_H:\Lambda(m,n)\rightarrow W(m,n)$ such that
  \[D_H(f)=\sum^s _{i=1}f_iD_i, \] 
for all $f\in \Lambda(m,n  )$, where $f_i=\sigma (i')(-1) ^{\tau(i')\d(f)}D_{i'}(f)$ and $\tau(i)=\d(D_i)$ for all $i\in Y$. Note that $D_H$ is even and that
\[[ D_H(f),D_H(g)]=D_H(D_H(f)(g)),\]
for all $f,g\in \Lambda(m,n)$. Let
  \[H(m,n):=\Span_{\BF}\{D_H(f)\mid f\in \Lambda(m,n)\}.\]
Then $H(m,n)$ is an infinite-dimensional Lie superalgebra. Let 
  \[\overline{H}(m,n;\underline{t})= \Span_{\BF}\{D_H(f)\mid f\in \Lambda(m,n;\underline{t})\}.\]
Then the set
  \[H(m,n;\underline{t})=[\overline{H}(m,n;\underline{t}),\overline{H}(m,n;\underline{t})]\]
is a finite-dimensional simple Lie superalgebra, which is called the Hamiltonian superalgebra. Note that $H(m,n;\underline{t})$ possesses a $\BZ$-graded structure:
  \[H(m,n;\underline{t})=\bigoplus\limits_{r=-1}^{\xi-3}H(m,n;\underline{t})_r,\]
\\by letting $H(m,n;\underline{t})_r=\Span_{\BF}\{D_H(x^{(\alpha)}x^u)\mid|\alpha|+|u|=r+2\}$. For convenience, we use $H$ and $H_r$ denote $H(m,n;\underline{t})$ and its $\BZ$-graded subspace $H(m,n;\underline{t})_r$ respectively.

\section{The notion of super-biderivation}

In this section, we give the definition of skew-symmetric super-biderivations of Lie superalgebras and obtain some useful conclusions about the skew-symmetric super-biderivations.

Let $G$ be a Lie algebra over an algebraically closed field. Recall that a linear map $D : G \rightarrow G$ is a derivation of $G$ if
\[D([x,y])=[D(x),y]+[x,D(y)],\]
for all $x, y\in G$. And we call a bilinear map $\varphi : G \times G\rightarrow G$ is a biderivation of $G$ if the following axioms are satisfied:
\begin{align*}
	\varphi (x,[y,z])&=[\varphi (x,y),z]+[y,\varphi (x,z)],\\
	\varphi([x,y],z)&=[\varphi(x,z),y]+[x,\varphi(y,z)],
\end{align*} 
for all $x, y, z\in G$. A biderivation $\varphi$ is called skew-symmetric if $\varphi(x, y) = -\varphi(y, x)$ for all $x, y \in G$.

A Lie superalgebra is a vector superspace $L=L_{\bar{0}}\,\oplus\, L_{\bar{1}}$ with an even bilinear mapping $[\cdot,\cdot]:L\times L\rightarrow L$ satisfying the following axioms:
\begin{align*}
	[x,y]&=-(-1)^{\d(x)\d(y)}[y,x],\\
	[x,[y,z]]&=[[x,y],z]+(-1)^{\d(x)\d(y)}[y,[x,z]],
\end{align*}
for all $x,y,z\in L$. Recall that a linear map $D:L\rightarrow L$ a superderivation of $L$ if 
\[D([x,y])=[D(x),y]+(-1)^{\d(D)\d(x)}[x,D(y)],\]
for all $x,y\in L$. Meanwhile, we write $\Der_{\bar{0}}(L)$ (resp. $\Der_{\bar{1}}(L)$) for the set of all superderivations of $\BZ_2$-degree $\bar{0}$ (resp. $\bar{1}$) of $L$. Let $\Der(L)=\Der_{\bar{0}}(L)\oplus\Der_{\bar{1}}(L).$

A $\BZ_2$-homogeneous bilinear map $\phi$ of $\BZ_2$-degree $\gamma$ of $L$ is a bilinear map such that $\phi(L_\alpha,L_\beta)\subset L_{\alpha+\beta+\gamma}$ for any $\alpha,\beta\in\BZ_2$. Specially, we call $\phi$ is even if $\d(\phi)=\gamma=\bar{0}$.

\begin{Definition}\label{def3.1}
	We call a bilinear map $\phi : L\times L\rightarrow L$ is a skew-symmetric super-biderivation of L if the following axioms are satisfied:
	\begin{align}
		\phi(x,[y,z])&=[\phi(x,y),z]+(-1)^{(\d(\phi)+\d(x))\d(y)}[y,\phi(x,z)], \label{eq3.1}\\
		\phi(x,y)&=-(-1)^{\d(\phi)\d(x)+\d(\phi)\d(y)+\d(x)\d(y)}\phi(y,x),\nonumber
	\end{align}
	for all $x,y,z\in L$. 
\end{Definition}
\noindent Observe that if $\phi$ is a skew-symmetric super-biderivation, the equation below is obviously fulled. 
\begin{align}
	&\phi([x,y],z)=[x,\phi(y,z)]+(-1)^{(\d(\phi)+\d(z))\d(y)}[\phi(x,z),y]. \label{eq3.2}
\end{align}
\noindent Meanwhile, we write $\BDer_{\bar{0}}(L)$ (resp. $\BDer_{\bar{1}}(L)$) for the set of all skew-symmetric super-biderivations of $\BZ_2$-degree $\bar{0}$ (resp. $\bar{1}$) of $L$. Let $\BDer(L)=\BDer_{\bar{0}}(L)\oplus\BDer_{\bar{1}}(L).$

\begin{Lemma}\label{lem3.1}
	If the bilinear map $\phi_{\lambda}$ : $L\times L\rightarrow L$ is defined by
	\[\phi_{\lambda}( x,y ) =\lambda [ x,y],\]
	for all $x,y\in L$, where $\lambda \in \BF$, then $\phi_{\lambda}$ is a skew-symmetric super-biderivation of $L$. This class of super-biderivations is called inner. Write $\IBDer(L)$ for the set of all inner super-biderivations of $L$.
	
	\begin{proof}
		Due to $\phi_{\lambda}( x,y ) =\lambda [ x,y]$, we deduce that $\phi_{\lambda}$ is an even super-biderivation, i.e. $\d(\phi_\lambda)=\bar{0}$. By the skew-symmetry of Lie superalgebras, it is easy to see that 
		\[\phi_\lambda(x,y)=-(-1)^{\d(\phi_\lambda)\d(x)+\d(\phi_\lambda)\d(y)+\d(x)\d(y)}\phi_\lambda(y,x),\]
		for any $x, y \in L$. And it is easily obtain the following three equalities:
		\begin{align*}
			\phi_\lambda(x,[y,z])&=\lambda[x,[y,z]],\\
			[\phi_\lambda(x,y),z]&=\lambda[[x,y],z],\\
			(-1)^{(\d(\phi_\lambda)+\d(x))\d(y)}[y,\phi_\lambda(x,z)]&=(-1)^{\d(x)\d(y)}\lambda[y,[x,z]].
		\end{align*}
		Since the graded Jacobi identity $[x,[y,z]]=[[x,y],z]+(-1)^{\d(x)\d(y)}[y,[x,z]]$,  we have that
		\[\phi_\lambda(x,[y,z])=[\phi_\lambda(x,y),z]+(-1)^{(\d(\phi_\lambda)+\d(x))\d(y)}[y,\phi_\lambda(x,z)],\]
		for any $x,y,z \in L$. Similarly, it follows that     
		\[\phi_\lambda([x,y],z)=[x,\phi_\lambda(y,z)]+(-1)^{(\d(\phi_\lambda)+\d(z))\d(y)}[\phi_\lambda(x,z),y],\] 
		for any $x, y, z \in L$. The proof is completed.
	\end{proof}
	
\end{Lemma}

\begin{Lemma}\label{lem3.2}		
	Let $L$ be a Lie superalgebra. Suppose that $\phi$ is a skew-symmetric super-biderivation of $L$. Then
	\[[\phi(x,y),[u,v]]=(-1)^{\d(\phi)(\d(y))+\d(u)}[[ x,y] ,\phi( u,v)],\]
	for any $x,y,u,v\in L$.
	
	\begin{proof}
		Due to the Definition \ref{def3.1}, there are two different ways to compute $\phi([x, u], [y, v])$. From the equation \eqref{eq3.1}, it follows that
		\begin{align*}
			\phi ([x,u],[y,v])
			=\ &[\phi([x,u],y), v]+(-1)^{(\d(\phi)+\d(x)+\d(u))\d(y)}[y,\phi([x,u],v)] \\
			=\ &[[x,\phi(u,y)], v]+(-1)^{(\d(\phi)+\d(y))\d(u)}[[\phi(x, y), u], v] \\
			&+(-1)^{(\d(\phi)+\d(x)+\d(u))\d(y)}[ y,[ x,\phi (u, v)]] \\
			&+(-1)^{(\d(\phi)+\d(x)+\d(u))\d(y)+(\d(\phi)+\d(v))\d(u)}[y, [\phi(x, v), u]].
		\end{align*}
		According to the equation \eqref{eq3.2}, one gets		
		\begin{align*}
			\phi([x,u],[y,v])
			=\ &[x,\phi(u,[y,v])]+(-1)^{(\d(\phi)+\d(y)+\d(v))\d(u)}[\phi(x,[y,v]),u] \\
			=\ &[x,[\phi(u,y),v]]+(-1)^{(\d(\phi)+\d(u))\d(y)}[x,[y,\phi(u,v)]] \\
			&+(-1)^{(\d(\phi)+\d(y)+\d(v))\d(u)}[[\phi(x,y),v],u] \\
			&+(-1)^{(\d(\phi)+\d(y)+\d(v))\d(u)+(\d(\phi)+\d(x))\d(y)}[[y, \phi(x, v)], u].
		\end{align*}
		
		Comparing two sides of the above two equations, and using the graded Jacobi identity of Lie superalgebras, we have that
		\begin{equation}\label{eq3.3}
			\begin{aligned}
				&[\phi(x,y),[u,v]]-(-1)^{\d(\phi)(\d(y)+\d(u))}[[x,y],\phi(u,v)]\\
				=\ &(-1)^{\d(y)\d(u)+\d(y)\d(v)+\d(u)\d(v)}([\phi(x,v),[u,y]]-(-1)^{\d(\phi)(\d(v)+\d(u))}[[x,v],\phi(u, y)]).
			\end{aligned}		
		\end{equation}		
		And we set 
		\[f(x,y;u,v)=[\phi(x,y),[u,v]]-(-1)^{\d(\phi)(\d(y)+\d(u))}[[x,y],\phi(u,v)].\]
		According to the equation \eqref{eq3.3}, it easily seen that 
		\[f(x,y;u,v)=(-1)^{\d(y)\d(u)+\d(y)\d(v)+\d(u)\d(v)}f(x,v;u,y).\]		
		On the one hand, we have
		\begin{align*}
			f(x,y;u,v)&=-(-1)^{\d(u)\d(v)}f(x,y;v,u)\\
			&=-(-1)^{\d(u)\d(v)}(-1)^{\d(y)\d(v)+\d(y)\d(u)+\d(v)\d(u)}f(x,u;v,y)\\
			&= (-1)^{\d(y)\d(u)}f(x,u;y,v).
		\end{align*}
		On the other hand, we also get
		\begin{align*}
			f(x,y;u,v)&=(-1)^{\d(y)\d(u)+\d(y)\d(v)+\d(u)\d(v)}f(x,v;u,y)\\
			&=-(-1)^{\d(y)\d(u)+\d(y)\d(v)+\d(u)\d(v)}(-1)^{\d(u)\d(y)}f(x,v;y,u)\\
			&=-(-1)^{\d(u)\d(y)}f(x,u;y,v).
		\end{align*}		
		Hence, it follows at once that 
		\[f(x,y;u,v)=-f(x,y;u,v).\]		
		Due to the characteristic of $\mathbb{F}$ is different from 2, we obtain that $f(x,y;u,v)=0.$ So we have 
		\[[\phi(x,y),[u,v]]=(-1)^{\d(\phi)(\d(y)+\d(u))}[[x,y],\phi(u, v)].\]
		The proof is completed.
	\end{proof}
	
\end{Lemma}

\begin{Lemma}\label{lem3.3}		
	Let $L$ be a Lie superalgebra. Suppose that $\phi$ is a skew-symmetric super-biderivation of $L$. If $\d(x)+\d(y)=\bar{0}$, then
	\[[\phi ( x,y ) ,[ x,y]] =0,\]
	for any $x,y\in L$.
	
	\begin{proof}
		By Lemma \ref{lem3.2}, we can have
		\[[\phi(x,y),[x,y]]=(-1)^{\d(\phi)(\d(y)+\d(x))}[[x,y],\phi(x, y)].\]		
		In view of $\d(x)+\d(y)=\bar{0}$, it is following that
		\begin{align*}
			[\phi(x,y),[x,y]]&=[[x,y],\phi(x,y)]\\
			&=-[\phi(x,y),[x,y]].
		\end{align*}
		Therefore, we obtain $[\phi(x,y),[x, y] ]=0$.
	\end{proof}
	
\end{Lemma}

\begin{Lemma}\label{lem3.4}   	
	Let $L$ be a Lie superalgebra. Suppose that $\phi$ is a skew-symmetric super-biderivation of $L$. If $[x,y]=0$, then
	\[\phi(x,y)\in Z_L([L,L]).\]
	where $Z_L([L,L])$ is the centralizer of  $[L,L]$. 	
	
	\begin{proof}
		If $[x,y]=0$, then we obtain
		\begin{align*}
			[\phi(x,y),[u,v] ]&= (-1)^{\d(\phi)(\d(y)+\d(u))}[[x,y],\phi[u,v]]\\
			&= 0,
		\end{align*}
		for any $u, v \in L$. So we get that $\phi(x,y)\in Z_L([ L,L])$.
	\end{proof}
	
\end{Lemma}

\section{Skew-symmetry super-biderivation of $H(m,n;\underline{t})$}

In this section, we use the method of the weight space decomposition of $H$ with respect to the abelian subalgebra $T_H$ to prove that all skew-symmetric super-biderivations of $H$ are inner. For convenience, we use $A$ and $B$ denote $A(m,\underline{t})$ and $B(n)$.

Set $T_H=\Span_{\BF}\{D_H(x_ix_{i'})\mid i\in Y\}$. Obviously, $T_H$ is an abelian subalgebra of $H$. For any $D_H(x^{(\alpha)}x^u)\in H$, we have
  \[[D_H(x_ix_{i'}),D_H(x^{(\alpha)}x^u)]=\sigma(i)(\alpha_{i'}-\alpha_{i}+\delta_{(i'\in  u)}-\delta_{(i\in  u)})D_H(x^{(\alpha)}x^u),\]
where
  \[\delta_{(\text{P})}=
  \begin{cases}
  	\ 1 \qquad \text{P is ture},\\
  	\ 0 \qquad \text{P is false}.
  \end{cases}\]
Fixed $\alpha\in A$ and $u\in B$, we define a linear function $(\alpha+\langle u\rangle):T_H\rightarrow \BF$
  \[(\alpha+\langle u\rangle)(D_H(x_ix_{i'}))=\sigma(i)(\alpha_{i'}-\alpha_{i}+\delta_{(i'\in  u)}-\delta_{(i\in  u)}).\]
Further, $H$ have a weight space decomposition with respect to $T_H$:
\[H=\bigoplus_{(\alpha+\langle u\rangle)}H_{(\alpha+\langle u\rangle)},\]
where
  \begin{align*}
  	H_{(\alpha+\langle u\rangle)}=\Span_{\BF}\{D_H(x^{(\beta)}x^v)\in H \mid &[D_H(x_ix_{i'}),D_H(x^{(\beta)}x^v)]=\\
  	&\sigma(i)(\alpha_{i'}-\alpha_{i}+\delta_{(i'\in  u)}-\delta_{(i\in  u)})D_H(x^{(\beta)}x^v)\}.
  \end{align*}

\begin{Lemma}\label{lem4.1}  	
	Let $\phi$ be a skew-symmetric super-biderivation of $H$. If $[x,y] = 0$ for $x,y \in H$, we have 
	\[\phi(x,y) = 0.\]	 
	
	\begin{proof}
		Since $H$ is a simple Lie superalgebra, it is obvious that $H=[H,H]$ and $Z(H)=0$. By Lemma \ref{lem3.4}, if $[x,y]=0$ for $x, y \in H$, we obtain
		\[\phi(x,y)\in Z_H([H,H])=Z(H)=0.\]
		The proof is completed.
	\end{proof}
	
\end{Lemma}

\begin{Lemma}\label{lem4.2}		
	Let $\phi$ be a skew-symmetric super-biderivation of $H$. For $D_H(x_ix_{i'})\in T_H$ and $D_H(x^{(\alpha)}x^u)\in H$, we have
	\[\phi(D_H(x_ix_{i'}),D_H(x^{(\alpha)}x^u))\in H_{(\alpha+\langle u\rangle)}.\]	
	
	\begin{proof}
		By lemma \ref{lem4.1}, it follows that $\phi(D_H(x_ix_{i'}),D_H(x_jx_{j'}))=0$ for any $i,j\in Y$ from $[D_H(x_ix_{i'}),$  $D_H(x_jx_{j'})]=0$. Note that $\d(D_H(x_ix_{i'}))=\bar{0}$, then for any $D_H(x^{(\alpha)}x^u)\in H$, it is clear that		
		\begin{align*}
			&[D_H(x_lx_{l'}),\phi(D_H(x_ix_{i'}),D_H(x^{(\alpha)}x^u))]\\
			=&(-1)^{(\d(\phi)+\bar{0})\,\bar{0}}\big(\phi(D_H(x_ix_{i'}),[D_H(x_lx_{l'}),D_H(x^{(\alpha)}x^u)])\\
			&-[\phi(D_H(x_ix_{i'}),D_H(x_lx_{l'})),D_H(x^{(\alpha)}x^u)]\big)\\
			=& \sigma(l)(\alpha_{l'}-\alpha_{l}+\delta_{(l'\in  u)}-\delta_{(l\in  u)})\phi(D_H(x_ix_{i'}),D_H(x^{(\alpha)}x^u)).
		\end{align*}
		The proof is completed.
	\end{proof}
	
\end{Lemma}

\begin{Lemma}\label{lem4.3}
	All $\BZ_2$-homogeneous skew-symmetric super-biderivations of $H$ are even. 
	
	\begin{proof}
		Suppose $\phi$ is a $\BZ_2$-homogeneous skew-symmetric super-biderivation of $H$. Since any element of $H_{(\alpha+\langle u\rangle)}$ has the same $\BZ_2$-degree, $\phi(D_H(x_ix_{i'}),D_H(x^{(\alpha)}x^u))$ and $D_H(x^{(\alpha)}x^u)$ have the same $\BZ_2$-degree. As we know that $\d(D_H(x_ix_{i'}))=\bar{0}$, then $\phi$ is even.
	\end{proof}
    
\end{Lemma}

\begin{Corollary}\label{cor4.4}
	All skew-symmetric super-biderivations of $H$ are $\BZ_2$-homogeneous, and fuether even. 
\end{Corollary}

We want to prove that all skew-symmetric super-biderivations of $H$ are inner. First, we need to prove the conclusion works on some elements of $H$. Therefore we give some specific weight spaces about the weight space decompositions of $H$ with respect to $T_H$.

\begin{Lemma}\label{lem4.5}
	Let $i\in Y_0,j\in Y_1$. Then the following statements hold:
	\begin{align*}
		&(1)H_{(\varepsilon_i)}=\sum_{\alpha\in A,\bar{u}\in B}\BF D_H\big((\prod_{\substack{r\in Y_0\backslash \{i,i'\}\\ 
				\alpha _{r'}-{\alpha _r} \equiv 0 \pmod {p}}} x^{(\alpha_r\varepsilon_r)})x^{((\overline{\alpha_{i'}+1})\varepsilon_i)}x^{(\alpha_{i'}\varepsilon_{i'})}x^{\bar{u}}\big);\\
		&(2)H_{(\langle j \rangle)}=\sum_{\alpha\in A,\bar{u}\in B}\BF D_H\big((\prod_{\substack{r\in Y_0\\ 
				\alpha _{r'}-{\alpha _r} \equiv 0 \pmod{p}}} x^{(\alpha_r\varepsilon_r)})x_jx^{\bar{u}}\big);\\
		&(3)H_{(\varepsilon_i+\langle j \rangle)}=\sum_{\alpha\in A,\bar{u}\in B}\BF D_H\big((\prod_{\substack{r\in Y_0\backslash \{i,i'\}\\ 
				\alpha _{r'}-{\alpha _r} \equiv 0 \pmod{p}}} x^{(\alpha_r\varepsilon_r)})x^{((\overline{\alpha_{i'}+1})\varepsilon_i)}x^{(\alpha_{i'}\varepsilon_{i'})}x_jx^{\bar{u}}\big);&
	\end{align*}
    where $j$ and $j'$ are both in $\bar{u}$ for $j\in Y_1$, $\overline{q}$ enotes some integer and $\overline{q}\equiv q\pmod{p}$. 
\end{Lemma}

\begin{Lemma}\label{lem4.6}
	Let $\phi$ is a skew-symmetric super-biderivation of $H$. For any $i\in Y$, there is an element $\lambda_i\in\BF$ such that
	\[\phi(D_H(x_ix_{i'}),D_H(x_i))=\lambda_i[D_H(x_ix_{i'}),D_H(x_i)],\]
	where $\lambda_i$ is dependent on $i$.
	
	\begin{proof}
		For $i\in Y_0$, by Lemma \ref{lem4.5} (1),  we can suppose that
		\begin{align*}
			\phi(D_H(x_ix_{i'}),D_H(x_i))
			=\sum_{\alpha\in A,\bar{u}\in B}a(\alpha,\bar{u},i)D_H\big((\prod_{\substack{r\in Y_0\backslash \{i,i'\}\\ 
					\alpha _{r'}-{\alpha _r} \equiv 0 \pmod{p}}} x^{(\alpha_r\varepsilon_r)})x^{((\overline{\alpha_{i'}+1})\varepsilon_i)}x^{(\alpha_{i'}\varepsilon_{i'})}x^{\bar{u}}\big),
		\end{align*}				
	    where $a(\alpha,\bar{u},i)\in\BF$. By Lemma \ref{lem4.1}, for any $l\in Y\backslash\{i,i'\}$, we have
	    \begin{align*}
	    	0&=\phi(D_H(x_ix_{i'}),[D_H(x_l),D_H(x_i)])-[\phi(D_H(x_ix_{i'}),D_H(x_l)),D_H(x_i)]\\
	    	&=[D_H(x_l),\phi(D_H(x_ix_{i'}),D_H(x_i))]\\
	    	&=\sum_{\alpha\in A,\bar{u}\in B}a(\alpha,\bar{u},i)D_H\big(D_H(x_l)((\prod_{\substack{r\in Y_0\backslash \{i,i'\}\\ 
	    			\alpha _{r'}-{\alpha _r} \equiv 0 \pmod{p}}} x^{(\alpha_r\varepsilon_r)})x^{((\overline{\alpha_{i'}+1})\varepsilon_i)}x^{(\alpha_{i'}\varepsilon_{i'})}x^{\bar{u}})\big).
	    \end{align*}
        By computing the equation, we find that $a(\alpha,\bar{u},i)=0$ if $\alpha_{l'}>0$ or $l'\in\bar{u}$. Thus we can suppose that
        \begin{align*}
            \phi(D_H(x_ix_{i'}),D_H(x_i))
            =\sum_{\alpha\in A}a(\alpha,i)D_H(x^{((\overline{\alpha_{i'}+1})\varepsilon_i)}x^{(\alpha_{i'}\varepsilon_{i'})}).
        \end{align*}
        Since $\d(D_H(x_ix_{i'}))+\d(D_H(x_i))=0$, by Lemma \ref{lem3.3}, we have
        \begin{align*}
        	0&=[\phi(D_H(x_ix_{i'}),D_H(x_i)),[D_H(x_ix_{i'}),D_H(x_i)]]\\
        	&=-\sigma(i)[\phi(D_H(x_ix_{i'}),D_H(x_i)),D_H(x_i)]\\
        	&=\sigma(i)[D_H(x_i),\phi(D_H(x_ix_{i'}),D_H(x_i))]\\
        	&=\sum_{\alpha\in A}a(\alpha,i)\sigma(i)D_H(D_H(x_i)(x^{((\overline{\alpha_{i'}+1})\varepsilon_i)}x^{(\alpha_{i'}\varepsilon_{i'})})).
        \end{align*}
        By computing the equation, we find that $a(\alpha,i)=0$ if $\alpha_{i'}>0$. Thus we can suppose that
        \begin{align*}
        	\phi(D_H(x_ix_{i'}),D_H(x_i))
        	=\sum_{\alpha\in A}a(\alpha,i)D_H(x^{(\alpha^{\bar{1}}_i\varepsilon_i)}),
        \end{align*}
        where $\alpha^{\bar{q}}_l$ denotes some integer and $\alpha^{\bar{q}}_l\equiv q\pmod{p}$. Similarly, we can suppose that
        \begin{align*}
        	\phi(D_H(x_ix_{i'}),D_H(x_{i'}))
        	=\sum_{\beta\in A}a(\beta,i')D_H(x^{(\beta^{\bar{1}}_{i'}\varepsilon_{i'})}),
        \end{align*}
        Since
        \begin{align*}
        	0&=\phi(D_H(x_ix_{i'}),[D_H(x_i),D_H(x_{i'})])\\
        	&=[\phi(D_H(x_ix_{i'}),D_H(x_i)),D_H(x_{i'})]+[D_H(x_{i}),\phi(D_H(x_ix_{i'}),D_H(x_{i'}))],
        \end{align*}
        which is equivalent to that
         \begin{align*}
        	[D_H(x_{i'}),\phi(D_H(x_ix_{i'}),D_H(x_i))]=[D_H(x_{i}),\phi(D_H(x_ix_{i'}),D_H(x_{i'}))],
        \end{align*}
      we have
        \begin{align*}
        	\sum_{\alpha\in A}a(\alpha,i)\sigma(i')D_H(D_i(x^{(\alpha^{\bar{1}}_i\varepsilon_i)}))=\sum_{\beta\in A}a(\beta,i')\sigma(i)D_H(D_{i'}(x^{(\beta^{\bar{1}}_{i'}\varepsilon_{i'})})).
        \end{align*}
        By computing the equation, we find that $a(\alpha,i)=0$ if $\alpha^{\bar{1}}_i>1$. Thus we can suppose that
         \begin{align*}
        	\phi(D_H(x_ix_{i'}),D_H(x_i))
        	=a(i)D_H(x_i),
        \end{align*}
        Set $\lambda_i=-\frac{a(i)}{\sigma(i)}$. Since $[D_H(x_ix_{i'}),D_H(x_i)]=-\sigma(i)D_H(x_i)$,
        we conclude that
        \[\phi(D_H(x_ix_{i'}),D_H(x_i))=\lambda_i[D_H(x_ix_{i'}),D_H(x_i)],\]
        where $\lambda_i$ is dependent on $i$.
        
        For $j\in Y_1$, by Lemma \ref{lem4.5} (2),  we can suppose that
         \begin{align*}
        	\phi(D_H(x_jx_{j'}),D_H(x_j))
        	=\sum_{\alpha\in A,\bar{u}\in B}a(\alpha,\bar{u},j)D_H\big((\prod_{\substack{r\in Y_0\\ 
        			\alpha _{r'}-{\alpha _r} \equiv 0 \pmod{p}}} x^{(\alpha_r\varepsilon_r)})x_jx^{\bar{u}}\big),
        \end{align*}
        By Lemma \ref{lem4.1}, for any $l\in Y\backslash\{j,j'\}$, we have
        \begin{align*}
        	0&=\phi(D_H(x_jx_{j'}),[D_H(x_l),D_H(x_j)])-[\phi(D_H(x_jx_{j'}),D_H(x_l)),D_H(x_j)]\\
        	&=[D_H(x_l),\phi(D_H(x_jx_{j'}),D_H(x_j))]\\
        	&=\sum_{\alpha\in A,\bar{u}\in B}a(\alpha,\bar{u},j)D_H\big(D_H(x_l)((\prod_{\substack{r\in Y_0\\ 
        			\alpha _{r'}-{\alpha _r} \equiv 0 \pmod{p}}} x^{(\alpha_r\varepsilon_r)})x_jx^{\bar{u}})\big).
        \end{align*}
        By computing the equation, we find that $a(\alpha,\bar{u},i)=0$ if $\alpha_{l'}>0$ or $l'\in\bar{u}$. Thus we can suppose that
        \begin{align*}
        	\phi(D_H(x_jx_{j'}),D_H(x_j))
        	=a(j)D_H(x_j).
        \end{align*}
        Set $\lambda_j=-\frac{a(j)}{\sigma(j)}$. Since $[D_H(x_jx_{j'}),D_H(x_j)]=-\sigma(j)D_H(x_j)$,
        we conclude that
        \[\phi(D_H(x_jx_{j'}),D_H(x_j))=\lambda_j[D_H(x_jx_{j'}),D_H(x_j)],\]
        where $\lambda_j$ is dependent on $j$.
	\end{proof}
    
\end{Lemma}

\begin{Lemma}\label{lem4.7}
	Let $\phi$ is a skew-symmetric super-biderivation of $H$. For any $i\in Y_0,j\in Y_1$, there is an element $\lambda_{ij}\in\BF$ such that
	\[\phi(D_H(x_ix_{i'}),D_H(x_ix_j))=\lambda[D_H(x_ix_{i'}),D_H(x_ix_j)],\]
	where $\lambda$ depends on neither $i$ nor $j$.
	
	\begin{proof}
		By Lemma \ref{lem4.5} (3),  we can suppose that
		\begin{align*}
			\phi(D_H(x_ix_{i'}),D_H(x_ix_j))
			=\sum_{\alpha\in A,\bar{u}\in B}a(\alpha,\bar{u},i,j)D_H\big((\prod_{\substack{r\in Y_0\backslash \{i,i'\}\\ 
					\alpha _{r'}-{\alpha _r} \equiv 0 \pmod{p}}} x^{(\alpha_r\varepsilon_r)})x^{((\overline{\alpha_{i'}+1})\varepsilon_i)}x^{(\alpha_{i'}\varepsilon_{i'})}x_jx^{\bar{u}}\big),
		\end{align*}				
		where $a(\alpha,\bar{u},i,j)\in\BF$. By Lemma \ref{lem4.1}, for any $l\in Y\backslash\{i,i',j,j'\}$, we have
		\begin{align*}
			0&=\phi(D_H(x_ix_{i'}),[D_H(x_l),D_H(x_ix_j)])-[\phi(D_H(x_ix_{i'}),D_H(x_l)),D_H(x_ix_j)]\\
			&=[D_H(x_l),\phi(D_H(x_ix_{i'}),D_H(x_ix_j))]\\
			&=\sum_{\alpha\in A,\bar{u}\in B}a(\alpha,\bar{u},i,j)D_H\big(D_H(x_l)((\prod_{\substack{r\in Y_0\backslash \{i,i'\}\\ 
					\alpha _{r'}-{\alpha _r} \equiv 0 \pmod{p}}} x^{(\alpha_r\varepsilon_r)})x^{((\overline{\alpha_{i'}+1})\varepsilon_i)}x^{(\alpha_{i'}\varepsilon_{i'})}x_jx^{\bar{u}})\big).
		\end{align*}
		By computing the equation, we find that $a(\alpha,\bar{u},i,j)=0$ if $\alpha_{l'}>0$ or $l'\in\bar{u}$. Thus we can suppose that
		\begin{align*}
			\phi(D_H(x_ix_{i'}),D_H(x_ix_j))
			=\sum_{\alpha\in A}a(\alpha,i,j)D_H(x^{((\overline{\alpha_{i'}+1})\varepsilon_i)}x^{(\alpha_{i'}\varepsilon_{i'})}x_j).
		\end{align*}
		By Lemma \ref{lem3.2} and \ref{lem4.6}, we have
		\begin{align*}
			0&=[\phi(D_H(x_ix_{i'}),D_H(x_i)),[D_H(x_ix_{i'}),D_H(x_ix_j)]]\\
			&=[[D_H(x_ix_{i'}),D_H(x_i)],\phi(D_H(x_ix_{i'}),D_H(x_ix_j))]\\
			&=-\sigma(i)[D_H(x_i),\phi(D_H(x_ix_{i'}),D_H(x_ix_j))]\\
			&=-\sum_{\alpha\in A}a(\alpha,i,j)\sigma(i)D_H(D_H(x_i)(x^{((\overline{\alpha_{i'}+1})\varepsilon_i)}x^{(\alpha_{i'}\varepsilon_{i'})}x_j)).
		\end{align*}
		By computing the equation, we find that $a(\alpha,i,j)=0$ if $\alpha_{i'}>0$. Thus we can suppose that
		\begin{align*}
			\phi(D_H(x_ix_{i'}),D_H(x_ix_j))
			=\sum_{\alpha\in A}a(\alpha,i,j)D_H(x^{(\alpha^{\bar{1}}_i\varepsilon_i)}x_j),
		\end{align*}
	    By Lemma \ref{lem3.2} and \ref{lem4.6}, we have
	    \begin{align*}
	    	\lambda_{i'}\sigma(i)D_H(x_j)&=[\lambda_{i'}[D_H(x_ix_{i'}),D_H(x_{i'})],[D_H(x_ix_{i'}),D_H(x_ix_j)]]\\
	    	&=[\phi(D_H(x_ix_{i'}),D_H(x_{i'})),[D_H(x_ix_{i'}),D_H(x_ix_j)]]\\
	    	&=[[D_H(x_ix_{i'}),D_H(x_{i'})],\phi(D_H(x_ix_{i'}),D_H(x_ix_j))]\\
	    	&=\sigma(i)[D_H(x_{i'}),\phi(D_H(x_ix_{i'}),D_H(x_ix_j))]\\
	    	&=\sum_{\alpha\in A}a(\alpha,i,j)\sigma(i)D_H(D_H(x_{i'})(x^{(\alpha^{\bar{1}}_i\varepsilon_i)}x_j)).
	    \end{align*}		
		By computing the equation, we find that $a(\alpha,i,j)=0$ if $\alpha^{\bar{1}}_i>1$. Thus we can suppose that
		\begin{align*}
			\phi(D_H(x_ix_{i'}),D_H(x_ix_j))
			=a(i,j)D_H(x_ix_j)=-\lambda_{i'}\sigma(i)D_H(x_ix_j),
		\end{align*}
		Since $[D_H(x_ix_{i'}),D_H(x_ix_j)]=-\sigma(i)D_H(x_ix_j)$. We conclude that
		\[\phi(D_H(x_ix_{i'}),D_H(x_ix_j))=\lambda_{i'}[D_H(x_ix_{i'}),D_H(x_ix_j)],\]
		By Lemma \ref{lem3.2} and \ref{lem4.6}, we have
		\begin{align*}			
			0=&[\phi(D_H(x_jx_{j'}),D_H(x_{j'})),[D_H(x_ix_{i'}),D_H(x_ix_j)]]\\
			    &-[[D_H(x_jx_{j'}),D_H(x_{j'})],\phi(D_H(x_ix_{i'}),D_H(x_ix_j))]\\
			=&(\lambda_{j'}-\lambda_{i'})[[D_H(x_jx_{j'}),D_H(x_{j'})],[D_H(x_ix_{i'}),D_H(x_ix_j)]]\\
			=&-(\lambda_{j'}-\lambda_{i'})\sigma(j)\sigma(i)[D_H(x_{j'}),D_H(x_ix_j)]\\
			=&(\lambda_{j'}-\lambda_{i'})\sigma(i)D_H(x_i).
		\end{align*}		
		Since $\sigma(i)D_H(x_i)\neq 0$, we have that $\lambda_{i'}=\lambda_{j'}$. From the arbitrariness of $i$ and $j$, we can claim that $\lambda_1=\lambda_2=\dots=\lambda_s$. Set $\lambda:=\lambda_1=\dots=\lambda_s$. Then we can conclude that for any $i\in Y_0,j\in Y_1$, there is an element $\lambda\in\BF$ such that
		\[\phi(D_H(x_ix_{i'}),D_H(x_ix_j))=\lambda[D_H(x_ix_{i'}),D_H(x_ix_j)],\]
		where $\lambda$ depends on neither $i$ nor $j$.		
	\end{proof}
	
\end{Lemma}

Then the conclusions of Lemma \ref{lem4.6} and \ref{lem4.7} can be rewritten as
\begin{align*}
	\phi(D_H(x_ix_{i'}),D_H(x_i))&=\lambda[D_H(x_ix_{i'}),D_H(x_i)],\\
	\phi(D_H(x_jx_{j'}),D_H(x_j))&=\lambda[D_H(x_jx_{j'}),D_H(x_j)],\\
	\phi(D_H(x_ix_{i'}),D_H(x_ix_j))&=\lambda[D_H(x_ix_{i'}),D_H(x_ix_j)],
\end{align*}
where $i\in Y_0,j\in Y_1$.

\begin{Theorem}
	Let $H$ be the Hamiltonian superalgebra $H(m,n;\underline{t})$ over the basic field $\BF$ of characteristic $p>2$. Then all skew-symmetric super-biderivations of $H$ are inner, i.e.
	\[\BDer(H)=\IBDer(H).\]
	\begin{proof}
		Suppose $\phi$ is a skew-symmetric super-biderivation of $H$. For any $D_H(f),D_H(g)\in H$, by Lemma \ref{lem3.2} and \ref{lem4.6}, for any $i\in Y$ we have
		\begin{align*}
			0=&[\phi(D_H(x_ix_{i'}),D_H(x_i)),[D_H(f),D_H(g)]]\\
			    &-[[D_H(x_ix_{i'}),D_H(x_i)],\phi(D_H(f),D_H(g))]\\
			 =&[\lambda[D_H(x_ix_{i'}),D_H(x_i)],[D_H(f),D_H(g)]]\\
			    &-[[D_H(x_ix_{i'}),D_H(x_i)],\phi(D_H(f),D_H(g))]\\
			 =&-\sigma(i)[D_H(x_i),\lambda[D_H(f),D_H(g)]-\phi(D_H(f),D_H(g))]\\
			 =&-(-1)^{\tau(i)}[D_{i'},\lambda[D_H(f),D_H(g)]-\phi(D_H(f),D_H(g))].
		\end{align*}
	    Since $Z_{H}(H_{-1})=\{E\in H\mid [E,D_i]=0,\forall\ i\in Y\}=H_{-1}$, we have that
	        \[\phi(D_H(f),D_H(g))-\lambda[D_H(f),D_H(g)]=\sum_{l\in Y}b_lD_l.\]
	    where $b_l\in \BF$. By Lemma \ref{lem3.2} and \ref{lem4.7} , we have
	    \begin{align*}
	    	0=&[\phi(D_H(x_ix_{i'}),D_H(x_ix_j)),[D_H(f),D_H(g)]]\\
	    	&-[[D_H(x_ix_{i'}),D_H(x_ix_j)],\phi(D_H(f),D_H(g))]\\
	    	=&[\lambda[D_H(x_ix_{i'}),D_H(x_ix_j)],[D_H(f),D_H(g)]]\\
	    	&-[[D_H(x_ix_{i'}),D_H(x_ix_j)],\phi(D_H(f),D_H(g))]\\
	    	=&\sigma(i)[D_H(x_ix_j),\sum_{l\in Y}b_lD_l]\\
	    	=&\sigma(i)(b_jD_H(x_i)-b_iD_H(x_j)).
	    \end{align*}         
         Since $\sigma(i)\neq 0,D_H(x_i)\neq0$ and $D_H(x_j)\neq0$, we have that $b_i=b_j=0$. Then $b_l=0,\ \forall\ l\in Y$. Hence, for any $D_H(f),D_H(g)\in H$, we have
             \[\phi(D_H(f),D_H(g))=\lambda[D_H(f),D_H(g)].\]
         Thus $\phi$ is an inner super-biderivation. From the arbitrariness of $\phi$, the proof is complete.
	\end{proof}
\end{Theorem}

%


\begin{thebibliography}{35}
	\setlength{\itemsep}{0mm} \small
	
	\bibitem{bai2023} W. Bai, W. Liu, Superbiderivations of Simple Modular Lie Superalgebras of Witt Type and Special Type, \emph{Algebra Colloq.} {\bf 30}(02) (2023), 181-192.
	
	\bibitem{Benkovic2009} D. Benkovi\v{c}, Biderivations of triangular algebras, \emph{Linear Algebra Appl.} {\bf 431}(9) (2009), 1587-1602.
	
	\bibitem{Bresar1993}M. Bre\v{s}ar, W. S. Martindale, C. R. Miers, Centralizing maps in prime rings with involution, J. Algebra {\bf 161}(2) (1993), 342-357.
	
	\bibitem{Bresar1995} M. Bre\v{s}ar, On generalized biderivations and related maps, \emph{J. Algebra} {\bf 172}(3) (1995), 764-786.
	
	\bibitem{Bresar2018} M. Bre\v{s}ar, K. Zhao, Biderivations and commuting linear maps on Lie algebras, \emph{J. Lie Theory} {\bf 28}(3) (2018), 885-900.
	
	\bibitem{chang20191} Y. Chang, L. Chen, Biderivations and linear commuting maps on the restricted Cartan-type Lie algebras $W(n;\underline{1})$ and $S(n;\underline{1})$, \emph{Linear Multilinear Algebra} {\bf 67}(8) (2019), 1625-1636.
	
	\bibitem{chang20192} Y. Chang, L. Chen and X. Zhou, Biderivations and linear commuting maps on the restricted Cartan-type Lie algebras $H(n;\underline{1})$, \emph{Comm. Algebra} {\bf 47}(3) (2019), 1311-1326.
	
	\bibitem{chang2021} Y. Chang, L. Chen and Y. Cao, Super-biderivations of the generalized Witt Lie superalgebra $W(m, n;\underline{t})$, \emph{Linear Multilinear Algebra} {\bf 69}(2) (2021), 233-244.
	
	\bibitem{chen2016} Z. Chen, Biderivations and linear commuting maps on simple generalized Witt algebras over a field, \emph{Electron. J. Linear Algebra} {\bf 31}(1) (2016), 1-12.
	
	\bibitem{cheng2017} X. Cheng, M. Wang, J. Sun and H. Zhang, Biderivations and linear commuting maps on the Lie algebra, \emph{Linear Multilinear Algebra} {\bf 65}(12) (2017), 2483-2493.
	
	\bibitem{cheng2019} X. Cheng, J. Sun, Super-biderivations and Linear Commuting Maps on the Twisted $N = 2$ Superconformal Algebra, \emph{Internat. J. Algebra Comput} {\bf 29}(07) (2019), 1235-1247.
	
	\bibitem{Dilxat2023} M. Dilxat, S. Gao and D. Liu. Super-biderivations and Post-Lie Superalgebras on Some Lie Superalgebras, \emph{Acta. Math. Sin. English Ser.} {\bf 39} (2023), 1736-1754. 
	
	\bibitem{du2013} Y. Du, W. Yu, Biderivations of generalized matrix algebras, \emph{Linear Algebra Appl.} {\bf438}(11) (2013), 4483-4499.
	
	\bibitem{fan2017} G. Fan, X. Dai, Super-biderivations of Lie superalgebras, \emph{Linear Multilinear algebra} {\bf 65}(1) (2017), 58-66.
	
	\bibitem{Ghosseiri2013}N.M. Ghosseiri, On biderivations of upper triangular matrix rings, \emph{Linear Algebra Appl.} {\bf 438}(1) (2013), 250-260.
	
	\bibitem{han2016} X. Han, D. Wang and C. Xia, Linear commuting maps and biderivations on the Lie algebras $W (a, b)$, \emph{J. Lie Theory} {\bf 26}(3) (2016), 777-786.
	
	\bibitem{Kac1977} V.G. Kac, Lie superalgebras, \emph{Adv. Math.} {\bf 26} (1977), 8-96.
	
	\bibitem{li2018} W. Li, X. Tang and J. Yuan, Super-biderivations and linear super-commuting maps on the super W-algebra $\widetilde {W}(2, 2) $, \emph{Colloq. Math.} {\bf 153} (2018), 273-300.
	
	\bibitem{Maksa1987}G. Maksa, On the trace of symmetric bi-derivations. \emph{C.R. Math. Rep. Acad. Sci. Canada} {\bf 9}(6) (1987), 303-307.
	
	\bibitem{Scheunert1979} M. Scheunert, The theory of Lie superalgebras, \emph{Lecture Notes in Math.} {\bf 716} (1979), 3001-454.
	
	\bibitem{tang2017} X. Tang, Biderivations, linear commuting maps and commutative post-Lie algebra structures on W-algebras, \emph{Comm. Algebra} {\bf 45}(12) (2017), 5252-5261.
	
	\bibitem{tang20181}X. Tang, Biderivations and commutative post-Lie algebra structures on the Lie algebra $\mathcal{W}(a,b)$, \emph{Taiwanese J. Math.} {\bf 22}(6) (2018), 1347-1366.
	
	\bibitem{tang20182} X. Tang, Biderivations of finite-dimensional complex simple Lie algebras, \emph{Linear Multilinear algebra}, {\bf 66}(2) (2018), 250-259.
	
	\bibitem{tang20183} X. Tang,  X. Li, Biderivations of the twisted Heisenberg-Virasoro algebra and their applications, \emph{Comm. Algebra} {\bf 46}(6) (2018), 2346-2355.
	
	\bibitem{tang2020} L. Tang, L. Meng and L. Chen, Super-biderivations and linear super-commuting maps on the Lie superalgebras, \emph{Comm. Algebra} {\bf 48}(12) (2020), 5076-5085.
	
	\bibitem{Vukman1990}J. Vukman, Two results concerning symmetric biderivations on prime rings, \emph{Aequationes
	Math.} {\bf 40} (1990), 181-189.
	
	\bibitem{wang2013} D. Wang, X. Yu, Biderivations and linear commuting maps on the Schr$\ddot{\text{o}}$dinger-Virasoro Lie algebra, \emph{Comm. Algebra} {\bf 41}(6) (2013), 2166-2173.
	
	\bibitem{wang2011} D. Wang, X. Yu and Z. Chen, Biderivations of the parabolic subalgebras of simple Lie algebras, \emph{Comm. Algebra} {\bf 39}(11) (2011), 4097-4104.
	
	\bibitem{xia2016} C. Xia, D. Wang and X. Han, Linear super-commuting maps and super-biderivations on the super-Virasoro algebras, \emph{Comm. Algebra} {\bf 44}(12) (2016), 5342-5350.
	
	\bibitem{xu2015} H. Xu, L. Wang, The properties of biderivations on Heisenberg superalgebras, \emph{Mathematica Aeterna} {\bf 5}(2) (2015), 285-291.
	
	\bibitem{yuan2021} J. Yuan, L. Chen and Y. Cao, Super-Biderivations of Cartan Type Lie Superalgebras, \emph{Comm. Algebra} {\bf 49}(10) (2021), 4416-4426.
	
	\bibitem{yuan2018} J. Yuan, X. Tang, Super-biderivations of classical simple Lie superalgebras, \emph{Aequationes Math.} {\bf 92}(1) (2018), 91-109.
	
	\bibitem{zhang2006} J.H. Zhang, S. Feng, H.X. Li and R.H. Wu, Generalized biderivations of nest algebras, \emph{Linear Algebra Appl.} {\bf 418}(1) (2006), 225-233.
	
	\bibitem{zhao2020} X. Zhao, Y. Chang, X. Zhou and L. Chen, Super-biderivations of the contact Lie superalgebra $K(m, n;\underline{t})$, \emph{Comm. Algebra} {\bf 488}(8) (2020), 3237-3248.
	
	\bibitem{zhang2005}Y. Zhang, W. Liu, Modular Lie superalgebras. \emph{Science Press, Beijing}, 2005 (in Chinese).
	
\end{thebibliography}
\end{document}